\documentclass[11pt]{article}
\usepackage[utf8]{inputenc}
\usepackage{amsmath}
\usepackage{amsfonts}
\usepackage{amssymb}
\usepackage{graphicx}

\newcommand{\N}{\mathbb{N}}
\newcommand{\R}{\mathbb{R}}

\newcommand{\BE}{\begin{equation}}
\newcommand{\EE}{\end{equation}}
\newtheorem{propo}{{P}roposition}
\newtheorem{defi}{{D}efinition}
\newtheorem{lem}{{L}emma}
\newtheorem{theo}{{T}heorem}
\newtheorem{cor}{{C}orollary}

\author{Sylvain Arguill\`ere, Universit\'e Pierre et Marie Curie (Univ. Paris 6),\\ CNRS UMR 7598, Laboratoire Jacques-Louis Lions, F-75005, \\Paris, France ({\tt sylvain.arguillere@upmc.fr})}
\title{Approximation of sequences of symmetric matrices with the symmetric rank-one algorithm and applications}
\date{}

\begin{document}

\maketitle

\begin{abstract}
The symmetric rank-one update method is well-known in optimization for its applications in the quasi-Newton algorithm. In particular, Conn, Gould, and Toint proved in 1991 that the matrix sequence resulting from this method approximates the Hessian of the minimized function. Extending their idea, we prove that the symmetric rank-one update algorithm can be used to approximate any sequence of symmetric invertible matrices, thereby adding a variety of applications to more general problems, such as the computation of constrained geodesics in shape analysis imaging problems. We also provide numerical simulations for the method and some of these applications.
\end{abstract}

\section{Introduction}

Let $d$ be an integer, $f:\R^d\rightarrow\R$ be a $\mathcal{C}^2$ function, and consider the problem of minimizing $f$ over $\R^d$. A well-known efficient algorithm to numerically solve this minimization problem is Newton's method: starting at some point $x_0$, it considers the sequence
$$
x_{k+1}=x_k-h_kH(f)^{-1}_{x_k}\nabla f(x_k),
$$
with $\nabla f$ the gradient of $f$, $H(f)$ its Hessian, and $h_k>0$ some appropriate step. 

\medskip

However, very often the Hessian of $f$ is too difficult to compute, leading to the  introduction  of the so-called quasi-Newton methods. The method defines a sequence
$$
x_{k+1}=x_k-h_kB_k^{-1}\nabla f(x_k),
$$ 
where $(B_k)$ is a sequence of symmetric matrices such that
\BE
\label{bkN}
B_{k+1}(x_{k+1}-x_k)=\nabla f(x_{k+1})-\nabla f(x_k).
\EE
Indeed, since
$$
\begin{aligned}
\nabla f(x_{k+1})-\nabla f(x_k)&=\left(\int_0^1H(f)_{x_k+t(x_{k+1}-x_k)}dt\right)(x_{k+1}-x_k)\simeq H(f)_{x_{k}}(x_{k+1}-x_k),
                               \end{aligned}
$$
we get
$$
B_{k+1}(x_{k+1}-x_k)\simeq H(f)_{x_k}(x_{k+1}-x_k).
$$
It is then expected that $B_k$ is close to $H(f)_{x_k}$ in the direction $s_k=x_{k+1}-x_k$. See \cite{qn1,qn2,qn3} for more.

There are many ways to build a matrix sequence $(B_k)$ satisfying \eqref{bkN}. However, it was proved in \cite{CGT} and \cite{S} that some of these methods let $B_k$ approximate $H(f)_{x_k}$ in all directions instead of just one, i.e.
$$
\Vert B_k-H(f)_{x_k}\Vert\underset{k\rightarrow\infty}{\rightarrow }0\quad\Longrightarrow\quad
\Vert B_k-H(f)_{x_*}\Vert\Vert B_k-H(f)_{x_k}\Vert\underset{k\rightarrow\infty}{\rightarrow }0.
$$
with the additional assumption of the uniform linear independence of the sequence 
$$
s_k=x_{k+1}-x_k,
$$
a notion that will be recalled later. In \cite{CGT} for example, this is proved for the update of $B_k$ by
\BE
\label{bk}
y_k=\nabla f(x_{k+1})-\nabla f(x_k)=A_ks_k,\quad r_k=B_ks_k-y_k,\quad B_{k+1}=B_k+\frac{r_kr_k^T}{r_k^Ts_k},
\EE
with 
$$
A_k=\int_0^1H(f)_{x_{k}+t(x_{k+1}-x_k)}dt.
$$

In this paper, our aim is to generalize the approach in \cite{CGT} by defining the above symmetric rank-one algorithm for any sequences of symmetric matrices $(A_k)$ and vectors $(s_k)$, and derive a convergence result, opening a wider range of applications. 

\medskip

For instance, if a sequence $A_k$ converges to an invertible matrix $A_*$, we can use the above algorithm to approximate the inverse $A_*^{-1}$ of the limit $A_*$. Indeed, let $(e_0,\dots,e_{d-1})$ be the canonical vector basis of $\R^d$, and define the sequence $(s_k)$ in $\R^d$ by
\BE
\label{appli}
s_k:=A_ke_{k[d]},\quad y_k=A_k^{-1}s_k=e_{k[d]},
\EE
where $k[d]$ is the remainder of the Euclidean division of $k$ by $d$. This sequence is uniformly linearly independent, hence the sequence $B_k$ defined by \eqref{bk} will converge to $A_*^{-1}$. This convergence might be a little slow depending on the dimension $d$ and the rate of convergence of $A_k$, but $B_k$ is much easier to compute than $A_k^{-1}$.

\medskip

This can be used to compute geodesics constrained to embedded submanifolds of Riemannian spaces. Indeed, to obtain a geodesic between two fixed points of a submanifold, we need to find a converging sequence of maps $t\mapsto\lambda_k(t)$ given implicitly by an equation of the form (\cite{ATY})
$$
A_k(t)\lambda_k(t)=c_k(t),
$$
where $A_k(t)$ is a convergent sequence of symmetric, positive definite matrices of high dimension. The $\lambda_k$ are the Lagrange multipliers induced by the equations of the submanifold. It is very consuming to solve such a linear system for every time $t$ and every step $k$. Instead, we can take
$$
\lambda_k(t)=B_k(t)c_k(t),
$$
with $B_k(t)$ obtained by applying the symmetric rank-one algorithm described in the previous paragraph. This is particularly useful in Shape Spaces, where the studied manifolds have a very high dimension and a very complex metric. The present article was actually motivated by such a problem appearing in shape analysis, investigated in \cite{ATY}.

\medskip

This paper is structured as follows. We  give the general framework in Section \ref{sec1}, then state the main result after recalling two equivalent definitions of the uniform linear independence of a sequence of vectors in Section \ref{sec2}. Then, Section \ref{sec3} is dedicated to intermediary results that will, along with notions developed in Section \ref{sec4}, derive the proof of our theorem.

\section{Notations and symmetric rank-one algorithm}
\label{sec1}

Consider a sequence $(A_k)_{k\in\N}$ of real square symmetric matrices of size $d$. Assume that this sequence converges to some matrix $A_*$, i.e.
$$
\Vert A_k-A_*\Vert\underset{k\rightarrow\infty}{\rightarrow }0,
$$
where $\Vert\cdot\Vert$ is the operator norm on $M_d(\R)$ induced by the canonical Euclidean norm $\vert\cdot\vert$ on $\R^d$. Then define
$$
\eta_{k,l}=\sup_{k\leq i\leq l}\Vert A_i-A_k\Vert,\quad\text{and}\quad\eta_{k,*}=\sup_{ i\geq k}\Vert A_i-A_k\Vert
$$
for all $k\leq l\in\N$. Note that
$$
\forall k\leq l\in\N,\quad \eta_{k,l}\leq\eta_{k,*}\qquad\text{and}\qquad \eta_{k,*}\rightarrow 0\quad\text{as}\quad k\rightarrow\infty.
$$
Now let $(s_k)_{k\in\N}$ be a sequence of vectors of $\R^d$.

\medskip
\noindent
Recall that we want to find a somewhat simple sequence $(B_k)_{k\in\N}$ of symmetric matrices such that $B_k\rightarrow A_*$, using only $s_k$ and $y_k=A_ks_k$.

\medskip

We use the symmetric rank-one update method from \cite{CGT}. Start with $B_0=I_d$. Then, define for $k\in\N$
$$
y_k=A_ks_k,\quad r_k=(A_k-B_k)s_k=y_k-B_ks_k,
$$ 
then take
$$
B_{k+1}=B_k+\frac{r_kr_k^T}{r_k^Ts_k}.
$$
It is of course required that $r_k^T s_k\neq 0$ for every $k$.

\section{Main Result}

\label{sec2}

For every $k$, we have
$$
B_{k+1}s_k=B_ks_k+r_k=B_ks_k+y_k-B_ks_k=y_k,
$$
so
$$
A_{k}s_k=B_{k+1}s_k.
$$
The main idea is that if $A_k,$ $A_{k+1}$, $\dots$, $A_{k+m}$ are not too far apart (i.e., for $k$ large enough), we expect $B_{k+m}s_{k+i}$ to be relatively close to $A_{k+m}s_{k+i}$ for $i\leq m$. Then, if we can extract from every finite subsequence $(s_k,\dots,s_{k+m})$ a vector basis of $\R^d$, we will obtain the desired convergence.

\medskip

For a more precise statement, we next define the notion of uniform linear independence. The most intuitive and geometric definition is the following.

\begin{defi}
\label{uli1}
Take a sequence $s=(s_k)_{k\in\N}$ of vectors in $\R^d$, $d\in\N^*$, and let $m\geq d$ be an integer. Then $s$ is said to be $m$-uniformly linearly independent if for some constant $\alpha>0$, and for all $k\in\N$, there are $d$ integers $k\leq k_1<\dots<k_d\leq k+m$ such that
$$
\left\vert\det\left( s_{k_1},\dots,s_{k_d}\right)\right\vert\geq \alpha\vert s_{k_1}\vert\dots\vert s_{k_d}\vert.
$$
\end{defi}
In other words, from every finite segment of $(s_k)$ of length $m$, we can extract a linear basis ${s_{k_1}},\dots,{s_{k_d}}$ that will, once normalized, form a parallelepiped that does not become flat as $k$ goes to infinity.

\medskip

Another definition was given in \cite{CGT} after \cite{ORBOOK} as follows.
\begin{defi}
\label{uli2}
A sequence $s=(s_k)_{k\in\N}$ of vectors in $\R^d$, $d\in\N^*$, is said to be $(m,\beta)$-uniformly linearly independent, where $d\leq m\in\N$ and $\beta\in\R$, if for all $k\in\N$, there are $d$ integers $k\leq k_1<\dots<k_d\leq k+m$ such that
$$
\left\vert \lambda\left(\frac{s_{k_1}}{\vert s_{k_1}\vert},\dots,\frac{s_{k_d}}{\vert s_{k_d}\vert}\right)\right\vert\geq \beta,
$$
where $\lambda(M)$ is the complex eigenvalue of the square matrix $M$ with smallest module.
\end{defi}

\textbf{Remark:} A sequence $s=(s_k)$ in $\R^d$ is $(m,\beta)$-uniformly linearly independent for some $m\geq d$ and $\beta>0$ if and only it is $m$-uniformly linearly independent in the sense of Definition \ref{uli1}.  Indeed, let $v_1,\dots,v_d\in \R^d$, and denote 
$
V=\displaystyle\left(\frac{v_1}{\vert v_1\vert},\dots,\frac{v_d}{\vert v_d\vert}\right).
$
If $\vert\lambda(V)\vert\geq \beta>0,$ then
$
\det(V)\geq \beta^d,
$
 which proves the first part of the equivalence. On the other hand,
we know that the eigenvalue of $V$ with largest modulus has modulus less than
$
\displaystyle\sqrt{d}\max_{i=1,\dots,d}\frac{\vert s_{k_i}\vert}{\vert s_{k_i}\vert}=\sqrt{d}.
$
Now, assume that $\det(V)\geq \alpha>0$. Then
$
\vert\lambda(V)\vert\geq \frac{\alpha}{d^{\frac{d-1}{2}}},
$
ensuring the second part of the equivalence.

\medskip

This definition is sufficient to state our main result. 
\begin{theo}
Let $(A_k),\ (s_k),\ (y_k),\ (r_k)$ and $(B_k)$ be defined as in Section \ref{sec1}, with $(A_k)$ having a limit $A_*$. Assume that for some fixed constant $c>0$,
$$
\vert r_k^T s_k\vert\geq c\vert r_k\vert\vert s_k\vert.
$$
Then, for every $\beta>0$ such that $(s_k)$ is $(m,\beta)$-uniformly linearly independent in the sense of Definition \ref{uli2}, we have for all $k\in \N$ the quantitative estimates
\BE
\label{resss}
\Vert B_{k+m}-A_*\Vert\leq \left(1+\left(\frac{2+c}{c}\right)^{m+1}\right)\frac{\sqrt{d}}{\beta}\eta_{k,*}.
\EE
\end{theo}

The next sections are dedicated to the proof of this theorem.

\section{First estimates}
\label{sec3}

In this section, we give upper bounds on
$$
\left\vert(B_{k+m}-A_{k})\frac{s_k}{\vert s_k\vert}\right\vert,
$$ 
and deduce estimates on 
$$
\left\vert\frac{(B_{k+m}-A_{*})x}{\vert x\vert}\right\vert
$$
for a particular set of $x\in \R^d$.
\begin{propo}
\label{mainlemma}
Let $(A_k)_{k\in\N}$ be a sequence of real symmetric matrices in $M_{d}(\R)$, $d\in\N$, and $(s_k)$ be any sequence in $\R^d$. Define $y_k$, $B_k$ and $r_k$ as above. Assume that for some fixed constant $0<c\leq 1$ and for all $k\in\N$,
$$
r_k^T s_k\geq c\vert r_k\vert\vert s_k\vert.
$$
Then, for all $l\geq k+1$, 
$$
\vert (A_k-B_l)s_k\vert\leq \left(\frac{2+c}{c}\right)^{l-k-1}\eta_{k,l-1}\vert s_k\vert.
$$
\end{propo}

\textbf{Proof:} We prove this inequality by induction on $l$, with $k\in \N$ fixed. For $l=k+1$, we know that $B_{k+1}s_k=A_ks_k=y_k$, hence
$$
\vert (A_{k}-B_{k+1})s_k\vert=0.
$$
We will use the notation
$$
IH(l):=\left(\frac{2+c}{c}\right)^{l-k-1}\eta_{k,l-1}\vert s_k\vert,
$$ 
where $IH$ stands for Induction Hypothesis. Now, assume the result to be true for some $l\geq k+1$, i.e.
\BE
\label{IH}
\vert (A_k-B_l)s_k\vert\leq \left(\frac{2+c}{c}\right)^{l-k-1}\eta_{k,l-1}\vert s_k\vert=IH(l).
\EE
Let us prove that
$$
\vert (A_k-B_{l+1})s_k\vert\leq \left(\frac{2+c}{c}\right)^{l-k}\eta_{k,l}\vert s_k\vert=IH(l+1).
$$
Note that
\BE
\label{akbl}
\begin{aligned}
\vert (A_k-B_{l+1})s_k\vert&=\vert A_ks_k-(B_l+\frac{r_lr_l^T}{r_l^Ts_{l}})s_k\vert\\
&=\vert A_ks_k-B_ls_k-\frac{r_lr_l^Ts_k}{r_l^Ts_{l}}\vert\\
&\leq \vert (A_k-B_l)s_k\vert+\frac{\vert r_l\vert\vert r_l^Ts_k\vert}{c\vert r_l\vert\vert s_l\vert}\\
&\leq IH(l)+\frac{\vert r_l^Ts_k\vert}{c\vert s_l\vert}.
\end{aligned}
\EE
Let us find a bound for  $\frac{\vert r_l^Ts_k\vert}{c\vert s_l\vert}$, the second term of the right-hand side. First we have
$$
\begin{aligned}
\vert r_l^Ts_k\vert&=\vert y_l^Ts_k-s_l^TB_l s_k\vert\\ 
&\leq \vert y_l^Ts_k-s_l^Ty_k\vert+\vert s_l^T(y_k-B_ls_k) \vert\\
&=\vert y_l^Ts_k-s_l^Ty_k\vert+\vert s_l^T(A_k-B_l)s_k \vert\\ 
&\leq \vert y_l^Ts_k-s_l^Ty_k\vert+\vert s_l\vert IH(l).
\end{aligned}
$$
However, since $A_l$ is symmetric and $y_l=A_ls_l$,
$$
\vert y_l^Ts_k-s_l^Ty_k\vert=\vert s_l^T(A_l-A_k)s_k\vert\leq \eta_{k,l} \vert s_l\vert\vert s_k\vert,
$$
from which we deduce
$$
\vert r_l^Ts_k\vert\leq \eta_{k,l} \vert s_l\vert\vert s_k\vert+IH(l)\vert s_l\vert.
$$
Going back to Inequality \eqref{akbl}, we get
$$
\begin{aligned}
\vert (A_k-B_{l+1})s_k\vert&\leq IH(l)+\frac{\vert r_l^Ts_k\vert}{c\vert s_l\vert}\\ 
&\leq IH(l)+\frac{1}{c} \eta_{k,l} \vert s_k\vert+\frac{1}{c} IH(l)\\
&=(1+\frac{1}{c}) IH(l)+\frac{1}{c}\eta_{k,l}\vert s_k\vert\\
&=\frac{1+c}{c}\left(\frac{2+c}{c}\right)^{l-k-1}\eta_{k,l-1}\vert s_k\vert+\frac{1}{c}\eta_{k,l}\vert s_k\vert \\
&\leq \left(\frac{2+c}{c}\right)^{l-k}\eta_{k,l}\vert s_k\vert=HI(l+1),
\end{aligned}
$$
where the last inequality comes from the simple fact that $\eta_{k,l-1}\leq\eta_{k,l}.\ \square$

\medskip

This proposition shows that if $A_k,\ A_{k+1},\dots,\ A_l$ are not too far away from each other (i.e. if $\eta_{k,l}$ is small), then $B_ls_k$ stays quantifiably close to $A_{k}s_k$.

\medskip

Now, note that $\Vert A_*-A_k\Vert\leq \eta_{k,*}$, and 
$
\eta_{k,*}
$
decreases to $0$ as $k$ goes to infinity. Keeping the same assumptions, we obtain the following result.

\begin{cor}
\label{aprox}
Take $m,k\in \N$, and let $x\in \R^d$ be in the span of $s_k,\dots,s_{k+m}$. If
$$
\frac{x}{\vert x\vert}=\sum_{i=0}^m \lambda_i\frac{s_{k+i}}{\vert s_{k+i}\vert},\quad \lambda_0,\dots,\lambda_m\in\R,
$$
then
$$
\frac{\vert B_{k+m}x-A_*x\vert}{\vert x\vert}\leq \eta_{k,*}\left(1+\left(\frac{2+c}{c}\right)^{m+1}\right)\sum_0^m\vert\lambda_i\vert.
$$
\end{cor}
\textbf{Proof:} First, it follows from Lemma \ref{mainlemma} that
$$
\begin{aligned}
\frac{\vert B_{k+m}x-A_{*}x\vert}{\vert x\vert}
&\leq \sum_{i=0}^m \frac{\vert \lambda_{k+i}\vert}{\vert s_{k+i}\vert}\vert  B_{k+m}s_{k+i}-A_{*}s_{k+i}\vert\\
&\leq \sum_{i=0}^m \frac{\vert \lambda_{k+i}\vert}{\vert s_{k+i}\vert}\Big(\vert B_{k+m}s_{k+i}-A_{k+i}s_{k+i}\vert+\vert A_{k+m}s_{*}-A_{k+i}s_{k+i}\vert\Big)\\
&\leq \sum_{i=0}^m \vert \lambda_i\vert\left(\left(\frac{2+c}{c}\right)^{i+1}\eta_{k,k+m-1}+\eta_{k,*}\right).
\end{aligned}
$$

Now, if we take
$$
C(m)=\left(1+\left(\frac{2+c}{c}\right)^{m+1}\right)
$$
and use the fact that $\eta_{k,k+m}\leq \eta_{k,*}$, then we get
$$
\frac{\vert B_{k+m}x-A_*x\vert}{\vert x\vert}
\leq \eta_{k,*}C(m)\sum_0^m\vert\lambda_i\vert.
$$
The result follows. $\square$

\medskip

In particular, if we can let $k$ go to infinity while keeping
$
\displaystyle\sum_{i=0}^m \vert \lambda_i\vert
$
bounded, then we obtain $B_{k+m}x\rightarrow A_*x$. Thus, if we can do it for all $x\in \R^d$, we will have proved that $B_k\rightarrow A_*$.

\medskip

In other words, we need every normalized vector $x\in\R^d$ to be a uniformly bounded linear combination of $s_k,\dots,s_{k+m}$ as $k$ goes to infinity. In the next section of this paper, we will define a third notion of uniform linear independence of a sequence directly related to to this property and prove that it is equivalent to the previous definitions.

\section{Uniform $m$-span of a sequence and applications}

\label{sec4}

In order to investigate the subspace on which $B_k\rightarrow A_*$, we need a notion that is more precise  than uniform linear independence.

\begin{defi}
\label{usm}
Let $s=(s_k)_{k\geq 0}$ be a sequence in $\R^d$, and let $m\in \N$. We say that a vector $x$ in $\R^d$ is uniformly in the $m-$span of $s$ if for some fixed $\gamma_x>0$,
\BE
\label{unif}
\forall k\in\N,\quad \exists\lambda_0,\dots,\lambda_m\in\R\quad\frac{x}{\vert x\vert}=\sum_{i=0}^m \lambda_i\frac{s_{k+i}}{\vert s_{k+i}\vert}\quad \text{and}\quad \sum_{i=0}^m \vert \lambda_i\vert\leq \gamma_x.
\EE
We denote by $US_m(s)$ the set of all such vectors. 

$US_m(s)$ is a vector sub-space of $\R^n$. Moreover, there exists a constant $\gamma>0$ such that Property \eqref{unif} holds for all $x\in US_m(s)$ with $\gamma_x=\gamma$, i.e.
\BE
\label{gsm}
\exists\gamma>0,\quad\forall k\in\N,\ x\in US_m(s),\quad \exists\lambda_0,\dots,\lambda_m\in\R,\quad\frac{x}{\vert x\vert}=\sum_{i=0}^m \lambda_i\frac{s_{k+i}}{\vert s_{k+i}\vert}\quad \text{and}\quad \sum_{i=0}^m \vert \lambda_i\vert\leq \gamma.
\EE
\end{defi}

To prove the existence of $\gamma$ in \eqref{gsm}, it suffices to consider an orthonormal basis $(x_i)_i$ of $US_m(s)$, associated with some constants
 $(\gamma_{x_i})_{1\leq i\leq d},$ in Property \eqref{unif}. Then we can just take $\gamma=\gamma_{x_1}+\dots+\gamma_{x_d}$.
 
\medskip
\noindent
\textbf{Remark:} There holds 
$
\displaystyle US_m(s)\subset\bigcap_{k=0}^{\infty}\text{span}(s_k,\dots,s_{k+m}).
$

\medskip
\noindent
\textbf{Example:} Define the sequence $s=(s_k)$ by 
$$
s_k=\left\lbrace \begin{aligned}&e_{k[d]}&\text{when}&\quad k[d]\neq n-1,\\
&e_0+\frac{1}{k}e_{d-1}&\text{when}&\quad k[d]= n-1,
\end{aligned}\right.
$$ 
where $k[d]$ is the remainder of the Euclidean division of $k$ by $d$. Then
$$
US_m(s)=\left\lbrace\begin{aligned}&\lbrace 0\rbrace\quad \text{if}\quad 0\leq m\leq n-2\\
&\text{span}(e_0,\dots,e_{d-2})\ \text{otherwise}.
\end{aligned}\right. 
$$

Using this definition, a simple application of Corollary \ref{aprox} gives the following result.

\begin{propo}
\label{mprop}
Let $(A_k),\ (s_k),\ (y_k),\ (r_k)$ and $(B_k)$ be defined as in Section \ref{sec1}, assuming that $(A_k)$ has a limit $A_*$ and that $\vert r_k^T s_k\vert\geq c\vert r_k\vert\vert s_k\vert$ for some fixed constant $c>0$.

\noindent
Then, for every $m\in \N$
\BE
\label{res}
\sup_{x\in US_m(\gamma)}\frac{\vert B_{k+m}x-A_*x\vert}{\vert x\vert}\leq C(m)\gamma\eta_{k,*},
\EE
where $\gamma$ is taken from \eqref{gsm} and
$$
C(m)=\left(1+\left(\frac{2+c}{c}\right)^{m+1}\right).
$$
\end{propo}

Finally the main result follows by combining this proposition with the following lemma.

\begin{lem}
\label{uli3}
Let $s=(s_k)_{k\geq 0}$ be a sequence in $\R^d$, and let $m\in \N$. Then $s$ is $(m,\beta)$-uniformly linearly independent if and only if $US_m(s)=\R^d$. Moreover, we can take
$
\gamma=\frac{\sqrt{d}}{\beta}
$ 
in \eqref{gsm}.
\end{lem}
\textbf{Proof:} Let $v_1,\dots,v_d$ be linearly independent elements of $\R^d$ and define the invertible matrix
$$
V=\left(\frac{v_1}{\vert v_1\vert},\dots,\frac{v_d}{\vert v_d\vert}\right).
$$
Let $\Lambda=(\lambda_1,\dots,\lambda_d)\in \R^d$, and define $x\in\R^d$ a normalized vector such that
$$
x=\sum_{i=1}^d\lambda_i v_i=V\Lambda.
$$
Then
$$
\sum_{i=1}^d\vert\lambda_i\vert\leq\sqrt{d}\vert\Lambda\vert=\sqrt{d}\vert
V^{-1}x\vert\leq \frac{\sqrt{d}}{\vert\lambda(V)\vert}.
$$
This proves that if a sequence $s=(s_k)$ in $\R^d$ is $(m,\beta)$-uniformly linearly independent, then $US_m(s)=\R^d$ and we can take $\gamma_m(s)=\frac{\sqrt{d}}{\beta}$.

\medskip

On the other hand, take $x\in\R^d$ a normalized such that 
$$
V^{-1T}V^{-1}x=\frac{1}{\lambda(V)^2}x.
$$
Then, if denoting $(\lambda_1,\dots,\lambda_d)=\Lambda=V^{-1}x$,
$$
\frac{1}{\vert\lambda(V)\vert}=\vert\lambda(V)\vert\ \frac{1}{\vert\lambda(V)\vert^2}=\vert\lambda(V)\vert\vert V^{-1T}V^{-1}x\vert=\vert\lambda(V)\vert\vert V^{-1T}\Lambda\vert\leq \vert\Lambda\vert\leq \sum_{i=1}^d\vert\lambda_i\vert,
$$
which proves the converse. $\square$

\medskip

Our main result is proved.

\section{Examples of applications and numerical simulations}

In this section, after running numerical simulations  of the algorithm on random symmetric matrices, we check that the inverse of a sequence of matrices can indeed be approximated. Then we give an application for computing constrained geodesics between two fixed points in Riemannian manifolds.

All simulations were done using Matlab.

\subsection{Approximation of a sequence of matrices}

Here we test the algorithm on random symmetric matrices with coefficients generated by a normalized Gaussian law. Let $d\in \N^*$, which will denote the size of the matrices.

First we define a square symmetric matrix $A_*=\frac{1}{2}(M+M^T)$, where the entries of the $d\times d$ matrix $M$ were chosen at random using the normalized Gaussian law. Then, fix $0<\lambda<1$, and define the sequence $(A_k)$ of symmetric matrices by perturbating $A_*$ as follows
$$
A_k=A_*+\frac{\lambda^k}{2}(M_k+M_k^T),
$$
where $M_k$ is a matrix with random coefficients taken uniformly in $[0,1]$.

Obviously, $A_k\rightarrow A_*$ linearly as $k\rightarrow\infty$.

Now, we define the sequence $(B_k)$ thanks to the symmetric rank-one algorithm, starting with $B_0=I_d$, and the sequence $(s_k)$ by the formula
$$
s_k=e_{k\ \text{mod}\ d},\quad k\in n,
$$
where $(e_0,\dots,e_{d-1})$ is the canonical basis of $\R^d$.

\medskip

Using the classical norm given by the Euclidean product $\left<X,Y\right>=tr(X^TY)$ on the space of square matrices of size $n$, we give in the following table the distance between $B_k$ and $A_*$. We took $d=10,$ and several values of $\lambda$ for a various number of steps.

\bigskip

\centerline{\begin{tabular}{|c|c|c|c|c|}
\hline
\text{Number of steps}  & 10 & 20 & 50 & 100\\
\hline 
$\lambda=$0.9 & 4 & 2 & 0.1 &  0.005 \\ 
\hline 
$\lambda=$0.5  & 1 & 1e-3 & 1e-12 & 0 \\ 
\hline 
$\lambda=$0.1  & 0.02 & 2e-12 & 0 & 0 \\ 
\hline 
\end{tabular} }

\bigskip

In the case of the usual quasi-Newton method (\cite{ORBOOK}, the goal is to approximate the inverse of the Hessian $H(f)$ of some function $f$ on $\R^d$. For this, we get a sequence of points $x_k$ converging to the minimum $x_*$, and we define
$$
\begin{aligned}
s_k&:=x_{k+1}-x_k,\\ 
A_k&:=\int_{0}^1 H(f)_{x_k+ts_k}dt,\\
 y_j&:=A_ks_k=\nabla g(x_{k+1})-\nabla g(x_k).
 \end{aligned}
$$
Then $\lim_{k\rightarrow\infty}A_k= H(x_*).$
\smallskip

This is a perfect example of a situation where it is easy to compute $A_k s_k$ for some particular $s_k$, but where it might be much harder to compute the actual $A_k$ (or just $H(x_k)$ for that matter). A large number of numerical simulations showing the efficiency of the symmetric rank-one algorithm can be found in \cite{CGT}.

\subsection{Approximation of the inverse of a sequence} 
 
As mentioned in the introduction, another application is the computation of the inverse $A_*^{-1}$ of the limit $A_*$, provided  $A_*$ is invertible.

\smallskip

Indeed, consider the following sequences for the symmetric rank-one algorithm
\BE
\label{skm1}
s_k:=A_ke_{k\ \text{mod}\ d},\ y_k=A_k^{-1}s_k=e_{k\ \text{mod}\ d}.
\EE
Then the sequence $(s_j)$ is $(d,\beta)$-linearly independent for some $\beta>0$ (at least starting at some $k_0$ large enough). Therefore, the sequence $B_k$ will converge to $A_*^{-1}$. This convergence might be a bit slow depending on the dimension $d$ and the rate of convergence of $A_k$, but $B_k$ is much easier to compute than $A_k^{-1}$.

This can be useful when solving approximately converging sequences of linear equations, as we will show in the next section.

\medskip

In the following numerical simulation, we used the same sequence $(A_k)$ with random coefficients as in the previous section, with $(A_k)$ converging linearly to a random matrix $A_*$ with rate $\lambda=0.5$. We then computed the distance between $B_k$ and $A_*^{-1}$. We also added an extra test. Indeed, the form of the sequence $s_k$ in \eqref{skm1} has no reason to be particularly good (i.e. uniformly linearly independent with a nice constant). Therefore, we applied the algorithm by taking a random vector $y_k$ with coefficients taken along a normal Gaussian law at each step and $s_k=A_ky_k$.

\bigskip

\centerline{\begin{tabular}{|c|c|c|c|c|}
\hline
\text{Number of steps}  & 10 & 20 & 50 & 100\\
\hline
$y_k=e_{k\ \text{mod}\ d}$ & 0.5 (300) & 0.001 (0.2) & 1e-7 (1e-4) & 1e-13 (1e-11) \\ 
\hline 
$y_k$ random  & 1 (500) & 0.01 (0.5) & 1e-6 (1e-3) & 1e-13 (1e-10) \\ 
\hline 
\end{tabular}}

\bigskip

For each case, we performed twenty experiments. Each entry in the previous table gives the mean value of the distance between $B_k$ and $A_*^{-1}$, with the highest value obtained in parentheses. This number can be significantly larger than the mean because of the randomness of $A_*$, which can cause it to be almost singular, leading the algorithm to behave badly as the $s_k$ are less uniformly linearly independent with this method.

\medskip

This experiment shows that taking $y_k$ random is not as efficient as taking $y_k$ to periodically be equal to the canonical basis of $\R^d$.

\subsection{An application: constrained optimisation}

Consider the following control system on $\R^d$
$$
\dot{x}=K_{x(t)}u(t),\quad u\in \R^d,
$$
with $K_x$ a semi-positive symmetric matrix with coefficients of class $\mathcal{C}^2$. This corresponds to a sub-Riemannian control system where $u$ is the momentum of the trajectory and $K_x$ the co-metric at $x$ (\cite{MBOOK}). This is a very natural formulation for problems of shape analysis, see \cite{ATY}. 

\medskip

Take $C\in M_{l,d}(\R)$ such that the $l\times l$ matrix 
$$
A_x=CK_xC^T
$$ 
is invertible for every $x\in \R^d$, and take an initial point $x_0\in \R^d$ such that $Cx_0=0$.

\medskip

We consider the optimal control problem of minimizing 
$$
L(u)=\frac{1}{2}\int_0^1u(t)^TK_{x(t)}u(t)dt+g(x(1)),\quad \text{where}\quad \dot{x}=K_{x(t)}u(t),
$$
over all possible $u\in L^2([0,1],\R^d)$ such that
$$
CK_{x(t)}u(t)=0\quad a.e.\ t\in[0,1].
$$
This is the same as minimizing $L(u)$ over trajectories that stay in the sub-space $\ker(C)$. According to the Pontryagin Maximum Principle (\cite{ATY,TBOOK}), if $u$ is optimal for this constrained problem, then there exists $p\in W^{1,2}([0,1],\R^d)$ such that $p(1)+dg_{x(1}=0$, and
\BE
\label{geodeq}
\left\lbrace
\begin{aligned}
\dot{x}&=K_x\left(p-C'A_x^{-1}CK_xp\right),\\
\dot{p}&=-\frac{1}{2}\left(p-C'A_x^{-1}CK_xp\right)^TK_x\left(p-C'A_x^{-1}CK_xp\right).
\end{aligned}\right.
\EE
for almost every $t\in[0,1]$. Moreover, 
\BE
\label{geocost}
L(u)=\tilde{L}(p(0))=\frac{1}{2}\left(p(0)-C'A_{x(0)}^{-1}CK_{x(0)}p(0)\right)^TK_{x(0)}\left(p(0)-C'A_{x(0)}^{-1}CK_{x(0)}p(0)\right)+g(x(1)).
\EE
Since \eqref{geodeq} is an ordinary differential equation, and since the minimization of $L$ reduces to the minimization of $\tilde{L}$ with respect to the initial momentum $p_0=p(0)$. Then, the computation the gradient of $\tilde{L}$ requires solving an adjoint equation with coefficients depending on the derivatives of the right-hand side of \eqref{geodeq}. This is described in more details in \cite{ATY}.

\medskip

One of the most time-consuming aspects of this method is the computation, at each time step, of the inverse of $A_x$. Therefore, we applied the Quasi-Newton Algorithm as follows. 

For any $k\in \N$, define $y_k=e_{k\ \text{mod}\ d}$. We start with the initial momentum $p_0=0$, and let $B_0(t)=Id_l$ for all $t\in[0,1]$. Then, assuming we have constructed an initial momentum $p_k$ and a family of matrices $B_k(t),\ t\in[0,1]$, we use \eqref{geodeq} to compute a trajectory $x_k(t)$, replacing $A_{x}^{-1}$ by $B_k{t}$. Finally, at each time $t$, we define 
$$
\begin{aligned}
s_k(t)&=A_{x_k(t)}y_k,\\
r_k(t)&=B_k(t)s_k(t)-y_k,\\
B_{k+1}(t)&=B_k(t)+\frac{r_k(t)r_k^T(t)}{r_k^T(t)s_k(t)}.
\end{aligned}
$$
We can then compute the gradient of $\tilde{L}$ with an adjoint equation, where any derivative 
$$
\partial_x(A_{x_k(t)})^{-1}=-A_{x_k(t)}^{-1}\partial_xA_{x_k(t)}A_{x_k(t)}^{-1}
$$ 
is replaced by $-B_k(t)\partial_xA_{x_k(t)}B_k(t)$. This allows the minimization of $\tilde{L}$ using gradient descent or a regular quasi-Newton algorithm.

\medskip

As long as the algorithm gives a converging sequence of initial momenta $p_k$, the trajectories $x_k(t)$ will also converge to a trajectory $x_*(t)$, making each $A_{x_k(t)}$, with $t\in[0,1]$ fixed, a converging sequence, with invertible limit $A_*(t)$. Therefore, each $B_k(t)$, $t\in[0,1]$ fixed, converges to $A_*(t)$ as $k\rightarrow \infty$. In other words, as $k\rightarrow\infty$, we are indeed computing the true gradient of $\tilde{L}$.

\newpage

\bibliographystyle{plain}
\bibliography{bibli}

\begin{thebibliography}{1}

\bibitem{ATY}
S.~Arguill{\`e}re, E.~Tr{\'e}lat, A.~Trouv{\'e}, and L.~Youn{\`e}s.
\newblock Control theory and shape analysis.
\newblock {\em Preprint}, 2013.

\bibitem{CGT}
A.~R. Conn, N.~I. Gould, and P.~L. Toint.
\newblock Convergence of quasi-newton matrices generated by the symmetric rank
  one update.
\newblock {\em Mathematical Programming}, (2):177--195, 1991.

\bibitem{qn1}
J.~E. Dennis, Jr. and Jorge~J. Mor{\'e}.
\newblock Quasi-{N}ewton methods, motivation and theory.
\newblock {\em SIAM Rev.}, 19(1):46--89, 1977.

\bibitem{qn2}
J.~E. Dennis, Jr. and Robert~B. Schnabel.
\newblock {\em Numerical methods for unconstrained optimization and nonlinear
  equations}, volume~16 of {\em Classics in Applied Mathematics}.
\newblock Society for Industrial and Applied Mathematics (SIAM), Philadelphia,
  PA, 1996.
\newblock Corrected reprint of the 1983 original.

\bibitem{qn3}
Philip~E. Gill, Walter Murray, and Margaret~H. Wright.
\newblock {\em Practical optimization}.
\newblock Academic Press Inc. [Harcourt Brace Jovanovich Publishers], London,
  1981.

\bibitem{MBOOK}
R.~Montgomery.
\newblock {\em A tour of subriemannian geometries, their geodesics and
  applications}, volume~91 of {\em Mathematical Surveys and Monographs}.
\newblock American Mathematical Society, Providence, RI, 2002.

\bibitem{ORBOOK}
J.~M. Ortega and W.~C. Rheinboldt.
\newblock {\em Iterative solution of nonlinear equations in several variables},
  volume~30 of {\em Classics in Applied Mathematics}.
\newblock Society for Industrial and Applied Mathematics (SIAM), Philadelphia,
  PA, 2000.
\newblock Reprint of the 1970 original.

\bibitem{S}
G.~Schuller.
\newblock On the order of convergence of certain quasi-{N}ewton methods.
\newblock {\em Numer. Math.}, 23:181--192, 1974.

\bibitem{TBOOK}
Emmanuel Tr{\'e}lat.
\newblock {\em Contr\^ole optimal}.
\newblock Math\'ematiques Concr\`etes. [Concrete Mathematics]. Vuibert, Paris,
  2005.
\newblock Th{\'e}orie \& applications. [Theory and applications].

\end{thebibliography}

\end{document}